\documentclass[12pt]{article}
\usepackage{setspace,amssymb,latexsym,amsmath,amscd,epsfig,amsthm,wasysym,enumerate}

\DeclareMathOperator*{\un}{\bigcup}

\newcommand{\real}{\mathbb{R}}
\newcommand{\nat}{\mathbb{N}}

\newcommand{\thm}[1]{\vskip.25in\noindent \textbf{Theorem {#1}}\\}
\newcommand{\cor}[1]{\vskip.25in\noindent \textbf{Corollary {#1}}\\}
\newcommand{\lma}[1]{\vskip.25in\noindent \textbf{Lemma {#1}}\\}

\newcommand{\pr}[1]{\vskip.25in\noindent \textbf{Proposition {#1}}\\}
\newcommand{\de}[1]{\vskip.25in\noindent \textbf{Definition {#1}}\\}
\newcommand{\con}[1]{\vskip.25in\noindent \textbf{Conjecture {#1}}\\}

\newcommand{\cont}{ \begin{flushright}\includegraphics[scale=0.075]{ContradictionTurtle} \end{flushright}}

\title{Orientable Hamiltonian Embeddings of Hypercubes}
\date{}
\author{Richard Leyland}

\begin{document}
\maketitle
\begin{abstract}
A Hamiltonian embedding is an embedding of a graph $G$ such that the boundary of each face is a Hamiltonian cycle of $G$. It is shown that the hypercube graph $Q_n$ admits such an embedding on an orientable surface when $n$ is a power of 2. Basic necessary conditions on Hamiltonian embeddings for $Q_n$ and conjectures are made about other values of $n$.
\end{abstract}
\section{Introduction}
The hypercube graph $Q_n$ is an example of a graph with a large amount of symmetry. We shall define $Q_n$ in two different ways.
\de{1.1} The $n$-dimensional hypercube $Q_n$, for $n \geq 1$ is given by
\[ V(Q_n) = \{0,1\}^n \qquad E(Q_n) = \{xy: x \mbox{ and $y$ differ in exactly one coordinate} \}\] 

We also provide an alternate definition which will be the basis of how we work with the hypercube graph.
\de{1.2 (Cartesian Product of Graphs)}
Let $G$ and $H$ be graphs. The cartesian product of graphs $G \Box H$ is given by
\[V(G \Box H) = V(G) \times V(H), \]
\[E(G \Box H) = \{(x,y)(x',y'): \mbox{ either } x=x' \mbox{ and } yy'\in E(H) \mbox{ or } y=y' \mbox{ and } xx' \in E(G)\}. \]
To illustrate this concept and the interpretation of the hypercube graph that we will be using, we introduce the merge operator.
\de{1.3} Let $G$ and $G'$ be isomorphic graphs given by the isomorphism $f: G \to G'$. The \emph{merge} operator define a new graph $G \star G'$, which is given by
\[V(G\star G')= V(G) \sqcup V(G') \qquad E(G\star G') = E(G)\cup E(G') \cup \{xy: x\in G, y\in G', f(x)=y \}. \] 
We shall refer to the edges of $G \star G'$ in the subset $E(G) \cup E(G')$ as \emph{inside} edges and we shall refer to the subset $\{xy: x\in G, y\in G', f(x)=y \}$ as \emph{outside} edges.\\

Now we consider $G \Box K_2$. We begin by labelling $V(K_2)= \{1,2\}$. So, by the definition of the Cartesian Product,
\[(x,i)(y,j) \mbox{ iff } u=v \mbox{ and } i \neq j \mbox{ or } i=j \mbox{ and }uv\in V(G). \]
This can be seen to be equivalent to the merge operation when $G=G'$ and $f$ is the natural isomorphism. Now $G \Box H$ can be interpreted in a similar manner by considering the disjoint union of copies of $G$ indexed by vertices of $H$. That is if $V(H) = \{x_1, \dots, x_n \}$, then we take $G^{x_1},\dots, G^{x_n}$. For every edge $xy\in E(H)$, perform the merge operation on $G^x$ and $G^y$.\\

Under this interpretation we shall denote $e(G^x,G^y)$ to denote the outside edges joining $G^x$ to $G^y$. This interpretation of the Cartesian product gives more clarity to the alternate definition of the $n$-dimensional hypercube $Q_n$.

\de{1.4 (Alternate Definition of $Q_n$)} For $n=1$ define $Q_1=K_2$. Then define $Q_n$ recursively so that $Q_n = Q_{n-1} \Box K_2$.\\

It turns out that $G \Box H$ is always isomorphic to $H \Box G$, but they are not the same. So we are careful in insisting that $Q_n = Q_{n-1} \Box K_2$. From this definition we can see that:
\[ Q_n = \underbrace{K_2 \Box K_2 \Box \dots \Box K_2}_{n \mbox{ times}}\]
and
\[Q_{2n} = \underbrace{K_2 \Box \dots \Box K_2}_{2n \mbox{ times}}= \underbrace{K_2 \Box \dots \Box K_2}_{n \mbox{ times}} \Box \underbrace{K_2 \Box \dots \Box K_2}_{n \mbox{ times}}= Q_n \Box Q_n.\]

\pr{1.5} $Q_n$ is bipartite with equal sized partitions.
\begin{proof}
We define the weight of a binary vector, $w(x)$, to be the number of non zero entries of $x$. Since $x$ is a binary vector, $w(x)$ is given by $\sum_{i=1}^n x_i$. Now if $xy \in E(Q_n)$, then either $w(x) = w(y)+1$ or $w(x) = w(y)-1$. So, $w(x)$ and $w(y)$ have opposite parities. Therefore, we define a partition $P = \{x: w(x) \equiv 0 \pmod{2}\}$ and $Q = \{x: w(x) \equiv 1 \pmod{2} \}$. It can be seen that $G[P]$ and $G[Q]$ both contain no edges. Futhermore, $|P|=|Q|$.
\end{proof}

We shall refer to vertices of $Q_n$ as being odd or even using the bipartition above. If $M$ is a matching of $Q_n$ then each matching edge is incident to a vertex of odd weight and a vertex of even weight. We shall refer to these vertices as the \emph{odd} and \emph{even endpoints} of the edge respectively.

\pr{1.6} The hypercube graph $Q_n$ has the following properties:
\begin{enumerate}
\item $v(Q_n) = 2^n$
\item $Q_n$ is $n$-regular
\item $e(Q_n) = n2^{n-1}$
\end{enumerate}

\begin{proof}
We shall proceed by induction using the cartesian product definition of $Q_n$. All three properties are easily checked for $Q_1 = K_2$. Now we suppose that the result holds for $Q_{n-1}$. Since, $Q_n = Q_{n-1} \Box K_2$, it follows that the process doubles the number of vertices. So,
\[v(Q_n) = 2v(Q_{n-1}) = 22^{n-1} = 2^n. \]
Furthermore, in the cartesian product, each vertex receives exactly one additional neighbour. So as $Q_{n-1}$ is $n-1$ regular, it follows that $Q_n$ must be $n$ regular. Lastly, by the handshake lemma
\[n2^n = \sum_{x\in V(Q_n)} \deg(x) = 2e(Q_n) \Rightarrow e(Q_n) = n2^{n-1}.\]
\end{proof}

Now, due to the high amount of symmetry, $Q_n$ is Hamiltonian. In fact for sufficiently large $n$, $Q_n$ has a lot Hamiltonian cycles. It is a fun exercise to inductively show that $Q_n$ is hamiltonian. In 1954 Ringel showed that the edge set $Q_n$ can be partitioned into hamiltonian cycles if $n$ is a power of two. Alspach et al. showed that this holds for every $Q_n$ with $n>2$. It is also known that $Q_n$ admits a Hamiltonian Cycle Double Cover, that is a list of Hamiltonian cycles $C_1, \dots, C_n$ so that each edge is contained exactly two of the cycles. However, we shall be looking for an even stronger structure. We are looking for an embedding of $Q_n$ on an orientable surface so that each face is bounded by a Hamiltonian cycle. First, we revise the definitions of embeddings.\\

\section{Embeddings of Graphs}
In the study of planar graphs, we look at the boundaries of faces (if the graph is 2-connected, all the boundaries will be cycles). If $F_1, \dots F_n$ are the faces, then each edge is on the boundary of exactly two of the faces, simply because on either side of the edge are two faces. So one might think that a cycle double cover of a graph is sufficient for specifying a planar embedding. However, some graphs are simply to dense for planar embeddings. For this we turn to surfaces.\\

We shall restrict ourselves to locally euclidean surfaces. That is surfaces so that every point has a neighbourhood that is homeomorphic to $\real^2$. It is well known that there are two such kinds of surfaces: orientable and non-orientable. Orientable surfaces are those in which rotating clockwise is distinct from the opposite, while on non-orientable surfaces there is no such distinction. The orientable surfaces are determined by the number of "holes" they have (up to homeomorphism). So there is the sphere, torus, the two holed torus and so on. The number of "holes" of a surface is called the genus, though we can define it more rigourously as the minimum number of handles we need to add on to a sphere so that the resulting surface is homeomorphic to the given surface. Non-orientable surfaces are determined by the number holes, crosscaps\footnote{A crosscap is given by identifying part of the surface with the boundary of a Mobius Strip} they have, but we shall be dealing only with orientable surfaces in this paper.\\

Now, we formally define what en embedding actually is. A surface $S$ is defined to be a compact connected 2-manifold\footnote{An t-manifold is a t-dimensional locally Euclidean space}. For a graph $G = (V,E)$, an embedding is a representation on $S$ such that the vertices of $G$ are points on $S$ and the edges of $G$ are simple arcs whose endpoints are the vertices. We maintain that the endpoint of the arc corresponding to $e$, correspond to the endpoints of the edge $e$. No arc includes points associated with other vertices and no two arcs intersect at a point which is interior to either arc. Technically, we are dealing with equivalence classes of surfaces, but the structure above is preserved under homeomorphism, since we demand that the embedding is injective.\\

Now as $S$ is locally euclidean we can look at a neighbourhood about each vertex in an embedding. Say we're looking at a vertex $v$. We enumerate the edges incident to $v$, say $e_1, \dots ,e_t$. Each of these edges borders two faces. So we take one of the cycles that contains $e_1$ and see what other edge it contains, say $f$. Now there is exactly one other cycle that contains $f$. We see what other edge on the cycle is incident to $v$. Repeating this process gives us a permutation on the edges incident to $v$. We can do this for all vertices and get a collection of  permutations $\{ \pi_v: v\in V(G)\}$ called a \emph{rotation scheme}. We can see that on an orientable surface, each permutation needs to be a cycle. That is if $v$ is of degree $d$, then $\pi_v$ must be a $d$-cycle. Now, we quote a well-known result in topological graph theory.\\

\thm{2.1} Specifying an embedding of a graph $G = (V,E)$ on a surface $S$ is equivalent to specifying a rotation scheme. If the surface is orientable, then we have an embedding if and only if $\pi_v$ is a $\deg(v)$-cycle for all $v\in V(G)$.\\

Here we do not distinguish between a permutation and it's cyclic rotations. We see that specifying a rotation scheme is sufficient for concluding that $G$ has an embedding. To do this, we take an edge $e=xy$. Then $\pi_y(xy) = yz$. So we consider $\pi_z(yz)$ and so on. We are guaranteed to get a closed walk from this procedure. Furthermore, we get that each edge will be contained in two closed walks. To construct the embedding, we take a 2-cell for each cycle and "sew" them together. It takes a lot to show that if $\pi_v$ is a cycle for all $v$, then the resulting surface is orientable but the details can be found in [1]. In this paper we shall be looking for a particular type of embedding.

\de{2.2}  A \textbf{Hamiltonian Embedding} of a $2$-connected graph $G$ is an embedding $\Sigma$ such that the cycle bounding each face is a Hamiltonian Cycle of $G$.

\de{2.3} A \textbf{Orientable Hamiltonian Embedding} of a $2$-connected graph is a Hamiltonian embedding $\Sigma$ on an orientable surface.\\

We can see that this is already a lot of structure to impose on a graph. In order to construct a Hamiltonian embedding, we first need a cycle double cover $C_1, \dots C_t$ such that each $C_i$ is a Hamiltonian cycle. Furthermore, in order to make sure that the embedding is orientable, we need to verify that the rotation scheme is a collection of cycles. That is $\pi_v$ is a $\deg(v)$ cycle for all $v$. Now we shall show that for certain values of $n$, $Q_n$ admits an Orientable Hamiltonian Embedding.

\section{Hamiltonian embeddings of $Q_{2^n}$}

Now we get to the central result.

\thm{3.1} If $n$ is a power of $2$ then $Q_n$ admits a Hamiltonian embedding on an orientable surface. 

We shall do this by decomposing $Q_n$ into $n$ perfect matchings, say $M_1, M_2, \dots, M_n$. It is easy to see that the union of two perfect matchings is a collection of disjoint cycles as the resulting graph will be 2-regular. We will construct the matchings so that $M_i \cup M_{i+1}$ is a Hamiltonian cycle and $M_i \cap M_j = \emptyset$. This is the central focus of the paper.

\thm{3.2} For $Q_{2^n}$ there exists a collection of disjoint perfect matchings $M_1,\dots, M_{2^n}$ such that $M_i \cup M_{i+1}$ is a Hamiltonian cycle. Here the subscripts are taken modulo $2^n$.\\

Why does this suffice? It is clear the defining $C_i = M_i \cup M_{i+1}$, will give us a collection of Hamiltonian cycles and the each edge will be contained in exactly two cycles. However, we need to verify the rotation scheme. Fix a vertex v. Let $e_i$ be the unique edge in $M_i$ incident to $v$. Now $e_i \in C_i$ and $e_i \in C_{i-1}$. Traversing along each cycle, we see that the rotation scheme for each vertex will be $\pi_v = (e_1, \dots, e_{2^n})$ or $\pi_v = (e_1, e_{2^n}, \dots,e_2)$, up to cyclic rotation. That is $e_i$ must lie between $e_{i-1}$ and $e_{i+1}$. So we are guaranteed that each $\pi_v$ is a $2^n$-cycle and so we have a Hamiltonian embedding.\\

\thm{3.2}
For each $n\in\nat$ there is a set of $m=2^n$ disjoint perfect matchings $M_1, \dots, M_m$ of $Q_{m}$, such that $M_i \cup M_{i+1}$ is a Hamiltonian cycle. 
\begin{proof}
We shall proceed by induction on $n$. $Q_2$ is clearly planar as $Q_2 \cong C_4$ and clearly has the desired decomposition. Namely if $Q_2 =\{00,10,01,11\}$. So set $M_1 = \{\{00,01\},\{10,11\}\}$ and $M_2 =  \{\{00,10\},\{01,11\}\}$. Note that in the argument above, $C_1$ and $C_2$ are the same cycle. We will have to take special note of this case when we show the matching decomposition of $Q_4$ is disjoint.\\

Suppose for $m=2^n$, we have such a matching decomposition for $Q_m$, $M_1, \dots, M_m$. We wish to construct such a matching decomposition for $Q_{2m}$. We shall use the indexing notation above. That is $Q_m^x$, for $x\in v(Q_m)$ refers to the implicit subgraph of $Q_{2m}$. Now for each matching $M_i$, we shall refer to the corresponding matching of $Q_m^x$ by $M_i^x$. So we define $N_i = \un_{x\in V(Q_n)}M_i^x$. $N_i$ is a perfect matching as any $v\in V(Q_m)$ belongs to a unique $M_i^x$. Therefore by definition of $M_i^x$, there is a unique edge of $M_i^x$ incident to $v$. So there is a unique edge of $N_i$ incident to $v$.\\

Now we define another set of matchings $O_i$ using the following classification. Let $e_E(Q_m^x,Q_m^y)$ denote the edges of $e(Q_m^x,Q_m^y)$ whose endpoint in $Q_m^x$ is even. Similarly let $e_O(Q_m^x,Q_m^y)$ denote the edges of $e(Q_m^x,Q_m^y)$ whose endpoint in $Q_m^x$ is even. We see that each edge of $e_E(Q_m^x, Q_m^y)$ has an odd endpoint in $Q_m^y$ and each edge of $e_O(Q_m^x,Q_m^y)$ has an even endpoint in $Q_m^y$. Futhermore $e_E(Q_m^x,Q_m^y) = e_O(Q_m^x,Q_m^y)$\\

We define $C_i = M_i \cup M_{i+1}$. Let $C_i = (x_{i_1}, \dots, x_{i_m})$. We orient these cycles in a clockwise direction with respect to the 2-cell embedding. Since the resulting surface is orientable this is well defined and we have if $xy$ appears in $C_i$ and $C_j$, then $x$ appears before $y$ in $C_i$'s orientation, and $y$ appears before $x$ in $C_j$'s orientation (or the other way around). This must be the case since left and right do not change direction anywhere on the surface. Rotate each of these so that $x_{i_1} = x_{j_1}$ for all $i,j = 1,\dots, 2^n$.
We define \[O_i = \un_{j=1}^m e_E(Q_m^{x_{i_j}},Q_m^{x_{i_{j+1}}})\]

We claim $O_i$ is a perfect matching for all $i$. Fix a vertex $v \in V(Q_{2m})$. Each $v$ belongs to a unique $Q_m^x$. Now, there are two edges in $C_i$ incident to $x$. Say $wx$ and $xy$. That is $C_i = (\dots, w,x,y, \dots)$. By the definition of merge, we have that there is a unique edge in $e(Q_m^x,Q_m^y)$ incident to $v$. Now if $v$ is of even weight, this edge is in $e_E(Q_m^x,Q_m^y)$ and therefore in $O_i$.\\

By the definition of $O_i$, the only other set that could have an edge incident with $v$ is the set $e_E(Q_m^w,Q_m^x)$. But each edge in this set has an odd endpoint in $Q_m^x$, so none of these edges can be incident with v and therefore $O_i$ has no other edges incident with $v$. Similarly, if $v$ is of odd weight, then there is a unique edge in $e_E(Q_m^w,Q_m^x)$ incident to $v$ and there are no edges in $e_E(Q_m^x,Q_m^y)$ incident with $v$.\\

Now we shall show that $N_i \cup O_i$ is a collection of $m/2$ cycles of length $2^{m}$. For each $v\in V(Q_m)$, let $v^x$ be the corresponding vertex in $Q_m^x$. Let $C_i = (x_{i_1}, \dots, x_{i_m})$ and orient them as in the statement of the theorem. First consider a matching edge $vw$ of $M_i$ We shall consider the matching edges of $M_i^{x_{i_j}}$, which are of the form $v^{x_{i_j}}w^{x_{i_j}}$. Now there is a unique edge of $O_i$, connecting $v^{x_{i_j}}$ to either $v^{x_{i_{j+1}}}$ or $v^{x_{i_{j-1}}}$. Without loss, we assume the unique edge is connected to  $v^{x_{i_{j+1}}}$. Note that by the definition of $O_i$, the assumption for one value of $j$, implies that for all j $w^{x_{i_j}}$ is connected to $w^{x_{i_{j+1}}}$ by an edge of $O_i$\\

So traversing edges of $O_i$ and $N_i$ alternatively we get the following cycle
\[ (v^{x_{i_1}}v^{x_{i_2}}w^{x_{i_2}}w^{x_{i_3}},\dots, v^{x_{i_{2^m-1}}}v^{x_{i_{2^m-1}}}w^{x_{i_{2^m}}}) \]
As we only increase the lowest index every other edge, we see the cycle is of length $2^m$. Since we get a cycle for each edge in $M_i$, we get a collection of $m/2$ cycles in the union $N_i \cup O_i$.\\

Of course we wish to fix this. To do this we shall use the fact that $M_i \cup M_{i+1}$ is a Hamiltonian cycle of $Q_m$. So, $M_i^x \cup M_{i+1}^x$ is a Hamiltonian cycle of $Q_m^x$. So we now define $P_i = M_{i+1}^{x_{i_1}} \cup \un_{j=2}^m M_i^{x_{i_j}}$. We use the following lemma to show $O_i \cup P_i$ is a Hamiltonian cycle of $Q_{2m}$.

\lma{3.3} Let $C_1, \dots, C_n$ be disjoint cycles of $G$. Suppose there is a matching $M$ so that for each $e\in M$, there is a unique $i$ such that $e\in C_i$. Suppose there is another matching $M'$ such that $M \cup M'$ is a cycle, $M'$ is disjoint from $M$ and $M' \cap C_i = \emptyset$ for every $i$. Then $\un(C_i \backslash M)\cup M'$ is a cycle of $G$.
\begin{proof}
Let $e_i = M \cap C_i$. Fix an orientation on $M \cup M'$ and relabel $C_1, \dots, C_n$ so that the $e_i$ appear in order on the traversal. Let $e_i=x_iy_i.$ So $M \cup M' = \{x_1y_1x_2y_2,\dots \}$. Then $M' = \{y_1x_2,y_2x_3,\dots, y_nx_1 \}$. Fix orientation on each $C_i$ so that they start off with $x_i$ and end on $y_i$. I.e, $C_i = (x_i, \dots, y_i)$. Let $D_i = C_i \ \{y_ix_i\}$, the path from $x_i$ to $y_i$ using edges of $C_i$. We see by traversing
\[ x_1D_1y_1x_2D_2y_2 \dots x_nD_ny_nx_1\]
is a closed walk. The traversal is a cycle since each cycle is disjoint from other cycles. We see this cycle includes no edges of $M$, all edges of $M'$ and all edges of $C_i \backslash M$.
\end{proof}
Let $D_1,\dots, D_{\frac{m}{2}}$ be the disjoint cycles obtained from $N_i \cup O_i$. Set $M=M_i^{x_{i_1}}$ and $M' = M_{i+1}^{x_{i_1}}$ and apply the lemma.
Then
\begin{eqnarray*}
(N_i \cup O_i) \cup M_{i+1}^{x_{i_1}} \backslash M_i^{x_{i_1}} &=& O_i \cup ((N_i\backslash M_i^{x_{i_1}}\cup M_{i+1}^{x_{i_1}})) \\
&=& O_i \cup P_i
\end{eqnarray*}
Now we claim $P_i \cup O_{i+1}$ is also a Hamiltonian cycle. First we get $N_i \cup O_{i+1}$ is a collection of disjoint cycles by repeating the argument for $N_i \cup O_i$. Then using the lemma, we get that indeed $P_i \cup O_i$ is a Hamiltonian cycle.\\

Now all that remains is to show that the matchings are disjoint. By construction, $P_i \cap P_j = \emptyset$ and $P_i \cap O_j$ for all $i \neq j$.  All that remains to show is that $O_i \cap O_j = \emptyset$ for all $i \neq j$.\\

Now we note by induction, $C_i \cap C_j$ is non empty if and only if $i=j+1$ or $i=j-1$.\\

Now, our definitions of $O_i$ depend on $C_i$. So an edge belongs to $O_i$ if and only if it belongs to $e_{E}(Q_m^x,Q_m^y)$ where $xy\in C_i$ and i is even or belongs to $e_{O}(Q_m^x,Q_m^y)$ where $xy\in C_i$ and i is odd.\\ 

First we consider the case where $m=2$. We started off with $M_1$ and $M_2$. But we noted that $C_1 = C_2$. At this stage we'll have constructed \[O_1 = e_E(Q_2^{00},Q_2^{10})\cup  e_E(Q_2^{10},Q_2^{11})\cup  e_E(Q_2^{11},Q_2^{01})\cup  e_E(Q_2^{00},Q_2^{01})\]. But \[O_2=e_O(Q_2^{00},Q_2^{10})\cup  e_O(Q_2^{10},Q_2^{11})\cup  e_O(Q_2^{11},Q_2^{01})\cup  e_O(Q_2^{00},Q_2^{01})\]
Which can be readily seen to be disjoint.

We continue for $m\geq 4$.\\

We have the following chain of equivalences. The which follow from the definition of $O_i$.
\[ xy\in C_i\cap C_j\iff e\in e_{E} (Q_m^x,Q_m^y) \mbox{ and } e\in e_{E}(Q_m^y,Q_m^x) \iff  e\in O_i \cap O_j\]
We know by induction that $xy\in C_i\cap C_j$ only when $i=j+1$ or $i=j-1$. Futhermore, by the orientations specified we can say that without loss $x$ appears before $y$ in the orientation of $C_i$ and $y$ appears before $x$ in the orientation of $C_j$. So we have the edges are taken from $e_E(Q_m^x,Q_m^y)$ and $e_E(Q_m^y,Q_m^x) = e_O(Q_m^x,Q_m^y)$ respectively. Here we see the importance of the specification of orientations, otherwise we could have $xy$ appearing in that order in both orientations and  we'd have the set $e_E(Q_m^x,Q_m^y)$, in both matchings. But due to our choice of orientations we are dealing with $e_E(Q_m^x,Q_m^y)$ and $e_O(Q_m^x,Q_m^y)$, which are completely disjoint. So $O_i \cap O_j$ is empty.\\

\end{proof}

\section{Necessary Conditions for Hamiltonian Embeddings}
A quick condition that can be derived is that if we have $d$ faces, all of which are Hamiltonian Cycles, we must have that the graph is $d$-regular. For a stronger result we turn to Heawood's generalization of Euler's Formula
\[v-e+f=2-2g\]
Where $g$ is the smallest number where a graph $G$ admits an embedding on an orientable surface of genus $g$. So noting the above, we suppose $G$ is a $d$-regular graph of order $n$. Therefore if we have a Hamiltonian embedding,
\[v= n \qquad e=\frac{nd}{2} \qquad f=d.\]
By Heawood's formula,
\[g = 1 - \frac{nd-2n-2d}{4}\]
So we must have
\[nd \equiv 2(n+d)\pmod{4}\]
which is summarized by the following

\pr{4.1}
Let $G$ is a $d$-regular graph of order $n$. Then if $G$ admits an orientable Hamiltonian Embedding we must have 
\begin{enumerate}
\item At least one of $n$ or $d$ must be even.
\item If $n \not\equiv d \pmod{2}$, then either $n\equiv 2 \pmod{4}$ or $d \equiv 2 \pmod{4}$.
\item If both $n$ and $d$ are even, then $nd \not\equiv n+d \pmod{4}$
\end{enumerate} 

From this simple result, we get that
\cor{4.2}
$Q_n$ has a Hamiltonian embedding only if $n$ is even.\\

Here we discuss some conjectures we have made. 
\con{4.3}
If $Q_n$ has an (orientable) Hamiltonian embedding with facial boundaries $C_1,\dots, C_n$, then either $C_i \cap C_j$ is a perfect matching of $Q_n$ or $C_i \cap C_j = \emptyset$.\\

We can see that if we want the embedding to be orientable, then two cycles $C_i$ and $C_j$ cannot share incident edges. To see this, suppose $e=xy$ and $f=yz$ belong to $C_i$ and $C_j$. We get that $\pi_y$ contains the transposition $(ef)$, when written as a union of disjoint cycles. So the embedding cannot be orientable.\\

We shall deal with the \emph{weighted intersection graph} $W$. Where $V(W) = \{C_1,\dots, C_n\}$. An edge $C_iC_j$ is in $W$ if and only if $C_i$ and $C_j$ share at least one edge. The weight of $C_iC_j$ is the number of edges they have in their intersection. This idea gives rise to the following conjecture:

\con{4.4} If $Q_n$ has an orientable Hamiltonian embedding with facial boundaries $C_1,\dots, C_n$. Then the weighted intersection graph is a cycle with weights $2^{n-1}$ on all edges.\\

Furthermore, we have the strong conjecture:
\con{4.5}
$Q_n$ has a Hamiltonian embedding only if $n$ is a power of 2.

\section*{References}
\begin{enumerate}[{[1]}]
\item Bojan Mohar and Carsten Thomassen. \emph{Graphs on Surfaces}. Johns Hopkins Studies in the Mathematical Sciences. Johns Hopkins University Press, Baltimore, MD, 2001.
\item Douglas B. West. \emph{Introduction to Graph Theory}. Pearson, 2000.
\end{enumerate}
 \end{document}